\documentclass{article}
\font\tenmath=msbm10
\font\sevenmath=msbm7
\font\fivemath=msbm5
\newfam\mathfam \textfont\mathfam=\tenmath
\scriptfont\mathfam=\sevenmath \scriptscriptfont\mathfam=\fivemath

\usepackage{amsmath, moreverb}
\usepackage{amsfonts, psfrag}
\usepackage{amssymb, fancybox}
\usepackage[american]{babel}
\usepackage{graphicx}
\usepackage{flafter}
\usepackage[section]{placeins}

\newcommand{\cJ}{\mathcal{J}}
\newcommand{\cE}{\mathcal{E}}

\newcommand{\cX}{\mathcal{X}}

\newcommand{\cP}{\mathcal{P}}
\newcommand{\cB}{\mathcal{B}}

\newcommand{\field}[1]{\mathbb{#1}}

\newcommand{\N}{{\rm I}\kern-0.18em{\rm N}}

\newcommand{\R}{{\rm I}\kern-0.18em{\rm R}}
\newcommand{\h}{{\rm I}\kern-0.18em{\rm H}}
\newcommand{\K}{{\rm I}\kern-0.18em{\rm K}}
\newcommand{\p}{{\rm I}\kern-0.18em{\rm P}}
\newcommand{\E}{{\rm I}\kern-0.18em{\rm E}}
\newcommand{\Z}{{\rm Z}\kern-0.18em{\rm Z}}

\newcommand{\1}{{\rm 1}\kern-0.25em{\rm I}}

\newcommand{\X}{\field{X}}
\newcommand{\ud}{\mathrm{d}}

\newcommand{\pn}{\p_{\kern-0.25em n}}
\newcommand{\pnm}{\p_{\kern-0.25em n,m}}
\newcommand{\psubm}{\p_{\kern-0.25em m}}
\newcommand{\e}{\textrm{e}}
\newcommand{\symdiff}{%
  \mathbin{\text{\footnotesize$\bigtriangleup$}}}
\newcommand{\symdiffsmall}{%
  \mathbin{\text{\scriptsize $\bigtriangleup$}}}
\newcommand{\Leb}{\mathrm{Leb}_d}

\newtheorem{TH1}{Theorem}[section]

\newtheorem{prop}{Proposition}[section]

\newtheorem{defin}{Definition}[section]

\begin{document}

\title{Generalization error bounds in semi-supervised classification under the
  cluster assumption}

\author{{\sc Philippe Rigollet}  \thanks{Laboratoire de Probabilit\'es et Mod\`eles
al\'eatoires, UMR~7599, Universit\'e Paris 6, case 188, 4,
pl. Jussieu, F-75252 Paris Cedex 5, France. Email:
{\tt rigollet@ccr.jussieu.fr}. Part of this work was done when the
  author was visiting researcher at Department of Statistics,
       University of California,
       Berkeley, CA 94720-1776, USA (Fund by France-Berkeley fund).}
}

\date{\normalsize \today}

\maketitle

\begin{abstract}
We consider semi-supervised classification when part of the available data
is unlabeled. These unlabeled data can be useful for the classification
problem when we make an assumption relating the behavior of the regression
function to that of the marginal distribution. Seeger \cite{s00} proposed the
well-known \emph{cluster assumption} as a reasonable one. We propose a
mathematical formulation of this assumption and a method based on
density level sets estimation that takes advantage of it to achieve fast rates
of convergence both in the number of unlabeled examples and the number of
labeled examples.

\end{abstract}

\medskip


\section{Introduction}
\setcounter{equation}{0}
Semi-supervised classification has been of growing interest over the past few years
and many methods have been proposed. The methods try to give an answer
to the question: ``How to improve classification accuracy using unlabeled data
together with the labeled data?''. Unlabeled data can be used in
different ways depending on the assumptions on the model. There are two types
of assumptions. The first one consists in assuming that we
have a set of potential classifiers and we want to aggregate them. In
that case, unlabeled data is used to measure the \emph{compatibility} between
the classifiers and reduces the complexity of the resulting classifier
(see, e.g., \cite{bb05}, \cite{bm98}). The second approach is the one that we use 
here. It assumes that the data contains clusters that have
homogeneous labels and the unlabeled observations are used to identify these
clusters. This is the so-called \emph{cluster assumption}. This idea can be
put in practice in several ways giving rise to various methods. The
simplest is the one presented here: estimate the clusters, then
label each cluster uniformly. Most of these methods use Hartigan's \cite{h75} definition of clusters, namely
the connected components of the density level sets. However, they use a parametric (usually mixture) model to estimate the underlying density
which can be far from reality. Moreover, no generalization error bounds are
available for such methods. In the same spirit,  \cite{t99} and \cite{r00} propose methods
that learn a distance using unlabeled data in order to have intra-cluster
distances smaller than inter-clusters distances. The whole family of graph-based
 methods aims also at using unlabeled data to learn the distances
between points. The edges of the graphs reflect the proximity between
points. For a detailed survey on graph methods we refer to \cite{z05}. Finally, we
mention kernel methods, where unlabeled data are used to build the
kernel. Recalling that the kernel measures proximity between points, such
methods can also be viewed as learning a distance using unlabeled data
(see \cite{bch04}, \cite{cz05}, \cite{czs06}). 

The cluster assumption can be interpreted in another way, i.e., as the requirement that the
decision boundary has to lie in low density regions. This interpretation has been
widely used in learning since it can be used in the design of standard algorithms such
as Boosting \cite{aga02}, \cite{hbw04} or SVM \cite{bch04}, \cite{cz05}, which are closely related to kernel methods mentioned above. In these algorithms, a greater penalization is given to decision
boundaries that cross a cluster.
For more details, see, e.g., \cite{s00}, \cite{z05}, \cite{czs06}. 
Although most methods make, sometimes implicitly, the cluster assumption, no
formulation in probabilistic terms has been provided so far. The formulation that we
propose in this paper remains very close to its original text formulation and allows to
derive generalization error bounds. We also discuss what can and
cannot be done using unlabeled data. One of the conclusions is that
considering the whole excess-risk is too ambitious and we need to concentrate on
a smaller part of it to observe the improvement of semi-supervised
classification over standard classification.

\vspace{0.1in}
\noindent \emph{Outline of the paper.} After describing the model, we
formulate the cluster assumption and discuss why and how it can improve
classification performance in the next section. In Section~\ref{secpop}, we
study the population
case when the marginal density $p$ is known, to get an idea of our
target. Indeed, such a population case corresponds  in some way to the
case when the amount of unlabeled data is infinite. Section~\ref{main} contains the main result: we propose an algorithm
for which we derive rates of
convergence for the $\lambda$-thresholded excess-risk as a measure of
performance. An exemple of consistent density level set estimators is
given in Section~\ref{secDLSE}. Section~\ref{secdisc} is devoted to discussion on the
choice of $\lambda$  and possible improvements. Proofs of the
results are gathered in Section~\ref{proofs}.

\vspace{0.1in}
\noindent \emph{Notation.}
Throughout the paper, we denote by $c_j$ positive constants. We write
$\Gamma^c$ for the complement of the set $\Gamma$. For two
sequences $(u_p)_p$ and $(v_p)_p$
 (in that paper, $p$ will be $m$ or $n$), we write $u_p=O(v_p)$ if there
 exists a constant $C>0$ such that $u_p \le Cv_p$ and we write $u_p=\widetilde
 O(v_p)$ if $u_p \le C (\log p)^\alpha v_p$ for some constants $\alpha>0,
 C>0$. Thus, if $u_p=\widetilde
 O(v_p)$, we have $u_p=o(v_p p^\beta)$, for any $\beta>0$.

\section{The model}
\setcounter{equation}{0}
Let $(X,Y)$ be a random couple with joint distribution $P$, where $X
\in \cX \subset \R^d$ is a
vector of $d$ features and $Y \in \{0, 1\}$ is a label indicating the class to
which $X$ belongs. 
The distribution $P$ of the random couple $(X, Y)$ is completely determined by
the pair $(P_X, \eta)$ where $P_X$ is the marginal distribution of $X$ and
$\eta$ is the regression function of $Y$ on $X$, i.e., $\eta(x)\triangleq
P(Y=1 |X=x)$. The goal of classification is to predict the label $Y$
given the value of $X$, i.e., to construct a measurable function $g:\cX \to \{0, 1\}$
called a \emph{classifier}. The performance
of $g$ is measured by the average classification error
$$
R(g) \triangleq P\left( g(X) \neq Y \right)
$$
A minimizer of the risk $R(g)$ over all classifiers is given by the
\emph{Bayes classifier}  $g^\star(x)
=\1_{\{\eta(x) \ge 1/2\}}$, where $\1_{\{\cdot\}}$ denotes the indicator
function. Assume that we have a
sample of $n$ observations $(X_1, Y_1), \ldots, (X_n, Y_n)$ that are independent
copies of $(X,Y)$. An empirical classifier is a random function  $\hat g_n:
\cX \to \{0,1\}$ constructed on the basis of the
sample $(X_1, Y_1), \ldots, (X_n, Y_n)$. Since $g^\star$ is the best possible
classifier, we measure the performance of an empirical classifier $\hat g_n$
by its \emph{excess-risk}
$$
\cE (\hat g_n)= \E_n R(\hat g_n) - R(g^\star)\,,
$$
where $\E_n$ denotes the expectation with respect to the joint distribution of
$(X_1,
Y_1), \ldots, (X_n, Y_n)$. We denote hereafter by $\pn$ the corresponding probability.

In many applications, a large amount of unlabeled data is available as well
as a small set of labeled data $(X_1, Y_1), \ldots, (X_n, Y_n)$ and the goal of semi-supervised classification is
to use of the unlabeled data to improve the performance of 
classifiers. Thus, we observe two independent samples
$ \X_l = \left\{(X_1,Y_1), \ldots, (X_n, Y_n)\right\}$ and $\X_u=\left\{ X_{n+1},
  \ldots, X_{n+m}\right\}$, where $n$ is rather small and typically $m\gg n$. It is well
known that in order to make use of the additional
unlabeled observations,
we have to make an assumption on the dependence between the marginal
distribution of $X$ and the joint distribution of $(X,Y)$. Seeger  \cite{s00}
formulated the rather intuitive
 \emph{cluster assumption} as follows\footnote{the
  notation is adapted to the present framework}
\begin{quote}
Two points $x, x' \in \cX$ should have the same label $y$ if there is a path between them
which passes only through regions of relatively high $P_X$.
\end{quote}
This assumption, in its raw formulation cannot be exploited in the
probabilistic model since $(i)$ the labels are random variables $Y, Y'$ so that the
expression ``should have the same label'' is meaningless unless $\eta$ takes
values in $\{0,1\}$ and $(ii)$ it is not clear what ``regions
of relatively high $P_X$'' are.
To match the probabilistic framework, we propose the following
 modifications
\begin{itemize}
\item[$(i)$] $P[ Y=Y' | X,X' {\rm connected}] \ge P[ Y \neq Y' | X,X' {\rm
    connected}] $, where ``connected'' means that there is the path between
    $X$ and $X'$ which passes only through regions of relatively high $P_X$.
\item[$(ii)$] Define ``regions of relatively high $P_X$'' in terms of
\emph{density level sets}.
\end{itemize}
We now need to precise the term \emph{ relatively high density}. Assume that $P_X$ admits a density $p$ with respect
to the Lebesgue measure on $\R^d$ denoted hereafter by $\Leb$. For a fixed $\lambda>0$, the
$\lambda$-level set of the density $p$ is defined by
\begin{equation}
\label{dls}
\Gamma(\lambda) \triangleq \left\{ x \in \cX: p(x) \ge \lambda\right\}\,.
\end{equation}

We are now in position to give a precise definition of the cluster
assumption. 
\begin{description}
\item[Cluster Assumption CA($\lambda$):] Fix $\lambda>0$ and  assume that the density
    level set $\Gamma=\Gamma(\lambda)$ has a countable number of connected
    components $T_j=T_j(\lambda),\ j=1,2, \ldots$.
Then the function $x \in \cX \mapsto \1{\{\eta(x) \ge
    1/2\}}$ takes a constant value on each of the $T_j, j=1,2, \ldots$.
\end{description}

Note that density level sets have the monotonicity property: 
$
\lambda \geq \lambda'$, implies $\Gamma( \lambda) \subset
\Gamma(\lambda')$. In terms of the cluster assumption, it means that when
$\lambda$ decreases to $0$, the assumption CA($\lambda$) becomes more restrictive. 
As a result, the parameter $\lambda$ can be considered as a level of confidence
characterizing to which extent the cluster assumption is valid for the
distribution $P$ and its choice is left to the user. For more details on the
choice of $\lambda$, see Section~\ref{secdisc}.

A question remains: what happens outside of the set $\Gamma(\lambda)$? Assume that we are in the problematic case,
  $P_X(\Gamma^c)=C>0$ such that the question makes sense. Since the cluster
  assumption says nothing about what happens outside of the set $\Gamma$, we
  can only perform supervised classification on $\Gamma^c$. Consider now a classifier $\hat g_{n,m}$ built from labeled and unlabeled
  samples $(\X_l, \X_u)$ pooled together. The excess-risk of $\hat g_{n,m}$ can be written
  (see \cite{dgl96})
$$
\cE(\hat g_{n,m})= \E_{n,m} \int_{\cX} |2 \eta(x) -1 | \1_{\{\hat g_{n,m}(x)\neq
  g^\star(x)\}}p(x) \ud x\,,
$$ 
where $\E_{n,m}$ denotes the expectation with respect to the pooled sample
$(\X_l, \X_u)$. We denote hereafter by $\pnm$ the corresponding probability.  Since, the unlabeled sample is of no help to
classify points in $\Gamma^c$, any reasonable classifier should be based on
the sample $\X_l$ so that $\hat
g_{n,m}(x)= \hat g_n (x), \ \forall\, x \in \Gamma^c$, and we have
\begin{equation}
\label{gammac}
\cE(\hat g_{n,m}) = \cE(\hat g_{n})\ge \E_n\int_{\Gamma^c}|2 \eta(x) -1 | \1_{\{\hat g_{n}(x)\neq
  g^\star(x)\}}p(x) \ud x\,.
\end{equation}
Since we assumed $P_X(\Gamma^c)=C>0$, the RHS of \eqref{gammac} is bounded
from below by the optimal rates of convergence that appear in supervised
classification. These rates are typically of the order $n^{-\alpha}, 1/2\le
\alpha \le 1$ (see e.g. \cite{mt99}, \cite{t04}, \cite{at05} and \cite{bbl05} for a comprehensive survey). Thus, unlabeled data do not improve the rate of convergence of this part of
the excess-risk. To observe the effect of unlabeled data on the rates of
convergence, we have to consider the $\lambda$-\emph{thresholded excess-risk}
of a classifier $\hat g_{n,m}$ defined by
\begin{equation}
\label{lamER}
\cE_{\lambda}(\hat g_{n,m}) \triangleq \E_{n,m} \int_{\Gamma(\lambda)} |2 \eta(x) -1 | \1_{\{\hat g_{n,m}(x)\neq
  g^\star(x)\}}p(x) \ud x\,.
\end{equation}
We will therefore focus on this measure of performance. Note that for
such a measure, we only need to consider classifiers $\hat g_{n,m}$ that are
defined on $\Gamma$. 

We now propose a method to obtain good upper bounds on this quantity, taking advantage of the cluster
assumption. The idea is to estimate the regions where the sign of $(\eta-1/2)$
is constant and make a majority vote on each region.

\vspace{0.2in}

\section{Results for known marginal distribution}
\label{secpop}
\setcounter{equation}{0}
Consider the ideal situation where the density $p$ is known and we observe only the labeled sample
$\X_l=\left\{(X_1, Y_1), \ldots, (X_n, Y_n)\right\}$. Fix $\lambda>0$ and
assume that  $\Gamma= \Gamma(\lambda)$ has a countable number of
connected components:
$$
\Gamma = \bigsqcup_{j\ge1} T_j\,,
$$
where the $T_j=T_j(\lambda)$ are non empty disjoint connected
sets. 
Under the cluster
assumption CA($\lambda$), the function $x \mapsto \eta(x) - 1/2$ has constant
sign on each $T_j$. Thus a simple and intuitive method for
classification is to perform a majority vote on each $T_j$. 

For any $j\ge 1$,
  define $\delta_j=\delta_j(\lambda) \ge 0$, $\delta_j \le 1$ by
$$
\delta_j\triangleq \int_{T_j}|2\eta(x) - 1| p(x) \ud x \,.
$$
The following assumption characterizes how far is $\eta$
from $1/2$ on every connected component $T_j$. 
\begin{description}
\item[Global Margin Assumption GMA($\lambda$):] There exists $\delta>0$ such
  that, for any $j\ge 1$, either $\delta_j=0$ or $\delta_j \ge \delta$.
\end{description}
Since $\sum_j \delta_j \le 1$, a direct consequence of the GMA is that only a
finite number of $\delta_j$ are positive. The GMA assumption imposes
that, on average over $T_j$, the regression function $\eta$ is away from
$1/2$ for any $j\ge 1$ such that $\delta_j>0$.  It describes the global behavior of $\eta$ on each connected component
$T_j$ as opposed to the standard margin assumption formulated in
\cite{mt99}  and \cite{t04} which we will call here \emph{local margin
  assumption} (LMA). Assumption LMA characterizes the local behavior of $\eta$ in a neighborhood of
$1/2$. In \cite{at05}, it is stated as follows 
\begin{description}
\item[Local Margin Assumption LMA:] There exist constants $C_0>0$ and
  $\alpha\ge 0$ such that
$$
P_X \left(0 < |2 \eta(X) -1| \le t \right) \le C_0 t^{\alpha}, \quad
\forall\, t \ge 0\,.
$$
\end{description}
It is straightforward that when there is only a finite number of connected
components $T_j, j=1, \ldots, J$ with non-zero Lebesgue measure,
GMA is a consequence of LMA. However we will see in
our analysis that the rates of convergence depend crucially on the value of
$\delta>0$, $j=1,2, \ldots$, while deriving GMA from LMA yields a $\delta$ depending on $C_0$. For this reason, it is natural to
introduce GMA instead of using the well known but less flexible
LMA.

We now define our classifier based on the sample $\X_l$ . For any $j\ge 1$, define the random variable
$$
Z_n^j \triangleq \sum_{i=1}^n \left(2Y_i - 1\right)\1_{\{X_i \in
  T_j\}}\,,
$$ 
and denote by $\hat g_n^j$  the function $\hat
g_n^j(x)=\1_{\{Z_n^j > 0\}}$ for all $x \in T_j$ .
Consider the classifier defined on $\Gamma$ by
$$
\hat g_n(x) = \sum_{j\ge 1} \hat g_n^j(x)\1_{\{x \in T_j\}}, \quad x \in \Gamma\,.
$$
The following theorem gives exponential rates of convergence for the
classifier $\hat g_n$ under CA($\lambda$).
\begin{TH1}
\label{thpop}
Fix $\lambda>0$ and assume that CA($\lambda$) holds. Then, the classifier
$\hat g_n$ satisfies
\begin{equation}
\label{ratepop}
\cE_{\lambda}(\hat g_n) \le 2 \sum_{j\ge 1} \delta_j \e^{-n\delta_j^2/2} \,.
\end{equation}
Moreover, if GMA($\lambda$) holds, inequality
\eqref{ratepop} reduces to
\begin{equation}
\label{ratepop2}
\cE_{\lambda}(\hat g_n) \le 2\e^{-n\delta^2/2}\,.
\end{equation}
\end{TH1} 
A rapid overview of the proof shows that the rate of convergence
$\e^{-n\delta^2/2}$ cannot be improved without further assumption. It will be
our target in semi-supervised classification. However, we need estimators of
the connected components $T_j, j\ge 1$. In the next section we provide
the main result on semi-supervised learning, that is when the density
$p$ is unknown but we can estimate it using the unlabeled sample $\X_u$. 

\section{Main result}
\label{main}
\setcounter{equation}{0}
We now deal with a more realistic case where the density $p$ is
 unknown and so are the density level sets which have to be estimated
 using the unlabeled sample $\X_u=\{X_1, \ldots, X_m\}$.  Fix $\lambda>0$ and
assume that  $\Gamma= \Gamma(\lambda)$ has a countable number of
connected components:
$$
\Gamma = \bigsqcup_{j\ge1} T_j\,,
$$
where the $T_j=T_j(\lambda)$ are non empty disjoint connected
sets. 
\subsection{Density level set estimation}
Assume that the density $p$ is uniformly bounded by a constant $L(p)$
 and that $\Leb(\cX)<\infty$, where $\Leb$ denotes the Lebesgue
 measure on $\R^d$. Denote by $\psubm$
 and $\E_m$ respectively the probability and the expectation w.r.t the sample
 $\X_u$ of size $m$.  Assume that
 for any $\lambda>0$, we use the sample $\X_u$ to construct an
 estimator $\hat G_m=\hat G_m(\lambda)$
 of $\Gamma=\Gamma(\lambda)$ satisfying
\begin{equation}
\label{estG}
\E_m\big[\Leb(\hat G_m \symdiff \Gamma) \big] \to 0, \quad  m \to + \infty.
\end{equation}
We call such estimators {\it consistent}
 estimators of $\Gamma$. However, the connected components of a
 consistent estimator of $\Gamma$ are not in general consistent
 estimators of the connected components of $\Gamma$. To ensure
 componentwise consistency, we have to make assumptions on the  connected
 component of $\Gamma$ and those of $\hat G$. 

\begin{figure}[h] \centering
\begin{minipage}[tr]{0.4\textwidth}
\psfrag{R}{\small $r_0$}
  \includegraphics[angle=0, width=\textwidth]{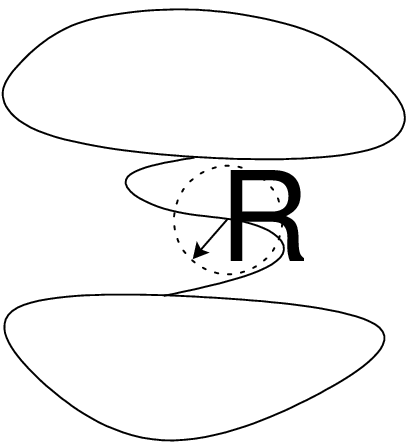}
\end{minipage}
\hspace{1cm}
\begin{minipage}[tl]{0.4\textwidth}
\psfrag{K}{\small $s_0/2$}
  \includegraphics[angle=0, width=\textwidth]{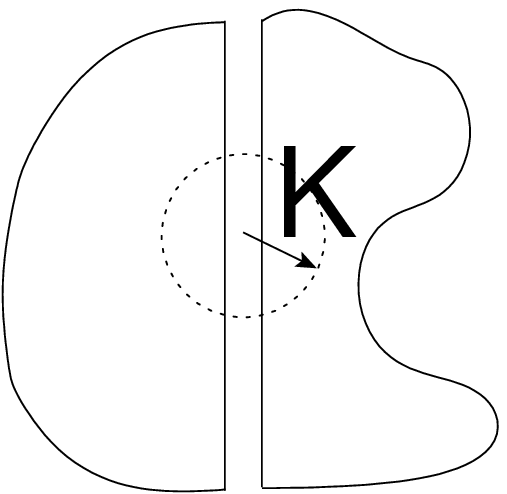}
\end{minipage}
\caption{Set that is 0-connected but not $r_0$-connected for any $r_0>0$ (left) and non-separated connected components (right).}\label{figsep}
\end{figure}

Let $\cB(x,
r)$ be the $d$-dimensional
closed ball of center $x
\in \R^d$ and radius $r>0$, defined by
$$
\cB(x, r)\triangleq \left\{ y \in \R^d: \|y-x\| \le r \right\}\,,
$$
where $\|\cdot \|$ denotes the Euclidean norm in $\R^d$. 
\begin{defin}
Fix $r_0 \ge 0$ and $c_0>0$. We say that a set $C \subset \R^d$
is $r_0$-\emph{connected} if for any $x, x' \in C$, there exists a continuous
map $f: [0,1] \to C$ such that $f(0)=x, f(1)=x'$and for any $t \in
[0,1]$ and any $r\le r_0$, we have
$$
\Leb\big(\cB(f(t),r)\cap C \big) \ge c_0r^d\,.
$$
A $0$-connected set is simply called connected or pathwise connected.
\end{defin}
This definition ensures that $\Gamma$ has no flat parts which allows to
exclude pathological cases such as the one presented on the left of Figure~\ref{figsep}. 
Now, define the distance $d_\infty$,
between two closed connected sets $C_1$ and $C_2$ by 
$$
d_\infty(C_1, C_2)= \min_{\substack{x \in C_1\\ y \in C_2}}\|x-y\|
$$
We say that a collection of connected sets $C_1, C_2, \ldots, $ is
 $s_0$-\emph{separated} if $d_\infty(C_j,
 C_{j'})\ge s_0, \ \forall j \neq j'$ for some $s_0 \ge 0$. If
 the connected components of $\Gamma$ are not $s_0$-separated for
 some $s_0>0$, cases such as the one presented on Figure~\ref{figsep} (right) could arise. In that case, two connected
 components and therefore two clusters are identified which is
 obviously not desirable. Therefore, the cluster assumption should not
 hold for that particular level $\lambda$ but it might hold for some
 $\lambda'\neq \lambda$.

Note that the performance of a density level
 estimator $\hat G_m$ is measured by the quantity
\begin{equation}
\label{symdiff}
\E_m\big[\Leb(\hat G_m \symdiff \Gamma) \big]=\E_m\big[\Leb(\hat G_m^c \cap \Gamma) \big]+\E_m\big[\Leb(\hat G_m \cap \Gamma^c) \big]\,.
\end{equation}
For some estimators, such as the penalized plug-in density level sets
estimators presented in Section~\ref{secDLSE}, we can prove that the dominant term in the RHS of
\eqref{symdiff} is $\E_m\big[\Leb(\hat G_m^c \cap \Gamma) \big]$. This
ensures that with high probability the estimator $\hat G_m$ is
included in $\Gamma$. We now give a precise definition of such
estimators.
\begin{defin}
Let $\hat G_m$ be an estimator of $\Gamma$ and fix $\alpha>0$. We say
that the estimator $\hat G_m$ is \emph{consistent from inside} at rate
$m^{-\alpha}$ if it satisfies
$$
\E_m\big[\Leb(\hat G_m \symdiff \Gamma) \big]=\widetilde O (m^{-\alpha})
$$ 
and 
$$
\E_m\big[\Leb(\hat G_m \cap \Gamma^c) \big]=\widetilde O (m^{-2\alpha})
$$
\end{defin}

For fixed $\alpha>0, \lambda>0$, let $\hat G_m \subset \cX$ be a
consistent from inside
estimator of $\Gamma$ at rate $m^{-\alpha}$. We begin by clipping $\hat G_m$ in
the following manner. Define the set 
$$
{\rm Clip}(\hat G_m)= \big\{ x \in \hat G_m: \Leb \big(\hat G_m \cap
\cB(x,(\log m)^{-1})\big) \le \frac{(\log m)^{-d}}{m^{\alpha}}\big\}.
$$
Note that $\Leb(\cX)<\infty$ yields
$$
\Leb\big({\rm Clip}(\hat G_m)\big) =\widetilde O(m^{-\alpha})
$$
and therefore the clipped set $\tilde G_m= \hat G_m \setminus {\rm
  Clip}(\hat G_m)$ is also consistent from inside at
rate~$m^{-\alpha}$. We now use only $\tilde G_m$. It is
straightforward that $\tilde G_m$
 can be decomposed into a finit number $\tilde J_m$ of connected
 components. We write for simplicity
\begin{equation}
\label{lab43}
\tilde G_m = \bigsqcup_{l\ge 1} \tilde T_l\,,
\end{equation}
where $\tilde T_l$ depends on $m$ and $\lambda$. Denote by
$\tilde H_k, k=1, 2, \ldots$, the family of sets such that
\begin{equation}
\label{lab44}
\bigsqcup_{l\ge1} \tilde T_l = \bigsqcup_{k\ge1} \tilde H_k\,,
\end{equation}
and $d_\infty (\tilde H_k,\tilde  H_{k'}) > 2(\log m)^{-1}$, $\forall \, k \neq k'$. It is not hard to see that the sets $\tilde H_k$ are uniquely
defined from $\tilde T_1, \tilde T_2, \ldots$. 
Let $\cJ$ be a subset of $\N^\star=\{1, 2, \dots\}$. Define $\kappa(j)=\{k :\tilde H_k
  \cap T_j \neq  \emptyset\}$ and let $D(\cJ)$ be the event on which the
  sets 
  $\kappa(j), j \in \cJ$ are reduced to 
  singletons $\{k(j)\}$  which are disjoint, i.e.,
\begin{equation}
\label{dkappa}
\begin{split}
D(\cJ) & \triangleq \Big\{ \kappa(j) =\{k(j)\}, k(j) \neq
k(j'),\    \forall\, j, j' \in \cJ,j \neq j' \Big\} \\
&=    \Big\{ \kappa(j) =\{k(j)\},  \ (T_j \cup \tilde H_{k(j)}) \cap (T_{j'}
  \cup \tilde H_{k(j')})= \emptyset,\    \forall\, j, j' \in \cJ,j \neq j' \Big\} \,.
\end{split}
\end{equation}
In other words, on the event $D(\cJ)$, there is a one-to-one
correspondence between the collection $\{T_j\}_{j \in \cJ}$ and
the collection $\big\{\{\tilde H_k\}_{k \in \kappa(j)}\big\}_{ j \in \cJ}$.
Componentwise convergence of $\tilde G_m$ to
$\Gamma$, is ensured when $D(\N^\star)$ has asymptotically overwhelming
probability. The following proposition gives an upper bound on the
probability of the complementary of $D(\cJ)$ under certain conditions
including the finiteness of $\cJ$.
\begin{prop}
\label{propcc}
Fix $r_0>0, s_0>0$ and let $\cJ$ be a subset of $\{1, 2, \ldots\}$
. Assume that $\{T_j\}_{j \in J}$ is a $s_0$ separated collection of
$r_0$-connected sets. Then, if $\hat
G_m$ is an  estimator of $\Gamma$ that is consistent from the inside
at rate $m^{-\alpha}$, we have
$$
\psubm \big(D^c(\cJ))= \widetilde{O}\big( m^{-\alpha}\big)\,.
$$
\end{prop}
The $r_0$-connectedness of all $T_j, j \in \cJ$ and $\Leb(\cX)<\infty$ entails that $\cJ$ is
necessarily finite. Nevertheless, the number of connected components
of $\Gamma$ can be infinite as long as there is only a finite number
of them for which $\delta_j=\int_{T_j} |2\eta-1|\ud P_X>0$.

 To estimate the homogeneous
regions, we will simply estimate the connected components of
$\Gamma$. In addition, when two connected components $T_j$ and $T_{j'}$ are close
with respect to the distance $d_\infty$, we merge\footnote{Merging two sets means here
  replacing them by their union} them into the same homogeneous region. 

It yields the following pseudo-algorithm.
\begin{center}
\label{r}
\ovalbox{
\begin{minipage}[h]{0.95\textwidth}
\centering  {\bf Pseudo-Algorithm}
\vspace{0.5cm}
\begin{enumerate}
\item Use the unlabeled data $\X_u$ to construct an estimator $\hat
  G_m$ of $\Gamma$ that is consistent from inside
at rate $m^{-\alpha}$. 
\item Define homogeneous regions as the unions of the connected
  components of $\tilde G_m=\hat G_m \setminus {\rm Clip}(\hat G_m)$ that are closer
  than $2(\log m)^{-1}$ for the distance $d_\infty$, accordinf to
  \eqref{lab43} and \eqref{lab44}.
\item Assign a single label to each estimated homogeneous region by majority
  vote on labeled data.
\end{enumerate}
\vspace{0.5cm}
\end{minipage}
}
\end{center}

This method translates into two distinct error terms, one term in $m$ and another
term in $n$. 
We apply our three-step procedure to build a  classifier $\tilde g_{n,m}$ based on the pooled
sample $(\X_l, \X_u)$. Fix $\lambda>0, \alpha> 0$ and let $\hat G_m$ be an
estimator of the density level set $\Gamma=\{p
\ge \lambda\}$, that is consistent from inside with rate
$m^{-\alpha}$.
For any
$k\ge 1$, define the random variable
$$
Z_{n,m}^k \triangleq \sum_{i=1}^n \left(2Y_i - 1\right)\1_{\{X_i \in
  \tilde  H_k \}}\,,
$$ 
where $\tilde H_k$ is defined in \eqref{lab44}. Denote by $\tilde g_{n,m}^k$  the function $\tilde
g_{n,m}^k(x)=\1_{\{Z_{n,m}^k > 0\}}$ for all $x \in \tilde H_k$ and consider the classifier defined on $\cX$ by
\begin{equation}
\label{classfin}
\tilde g_{n,m} \triangleq \sum_{k\ge1} \tilde g_{n,m}^k(x)\1_{\{x \in \tilde H_k\}}, \quad x \in \cX\,.
\end{equation} 
Note that the classifier $\tilde g_{n,m}$ assigns the label $0$ to any $x$
outside of $\tilde G_m$. This is a notational convention and we can assign
any value to $x$ on this
set since we are only interested in the $\lambda$-thresholded
excess-risk. Nevertheless, it is more appropriate to assign a label referring to  a
rejection, e.g., the values ``2''or  ``R'' (or any other value different from
$\{0,1\}$). The rejection meaning that this point should be classified
using labeled data only. However, when the amount of labeled data is too small,
it might be more reasonnable not to classify this point at all. This modification is of particular interest in the context of
classification with a rejection option when the cost of rejection is smaller
than the cost of misclassification (see, e.g., \cite{hw05}). 
\begin{TH1}
\label{mainth}
Fix $\lambda>0, \alpha>0, r_0>0$ and assume that CA($\lambda$)
holds. Consider an estimator $\hat G_m$ based on $\X_u$  that is consistent from inside with rate
$m^{-\alpha}$. Then if the connected components of $\Gamma(\lambda)$
are $r_0$-connected and $s_0$-separated, the classifier $\tilde  g_{n,m}$ defined in
\eqref{classfin} satisfies
\begin{equation}
\label{maineq}
\cE_{\lambda} \left(\tilde g_{n,m}\right) \le \widetilde O \left(\frac{m^{- \alpha}}{1-\theta}\right)  + \sum_{j\ge 1} \delta_j \e^{-n (\theta \delta_j)^2/2}\,,
\end{equation}
for any $0< \theta < 1$. Moreover, if GMA($\lambda$) holds, inequality
\eqref{maineq} reduces to
\begin{equation}
\label{maineq2}
\cE_{\lambda} \left(\tilde g_{n,m}\right) \le \widetilde O \left(\frac{m^{-
    \alpha}}{1-\theta}\right)  +  \e^{-n (\theta \delta)^2/2}\,.
\end{equation} 
\end{TH1}
Note that, since we often have $m \gg n$, the first term in the RHS of  \eqref{maineq} and \eqref{maineq2}  can be considered negligible so that we achieve an exponential rate of
convergence in $n$ which is almost the same (up to the constant
$\theta$ in the exponent) as in the case where the density $p$ is
known. The constant $\theta$ seems to be natural since it balances the
two terms.

\section{Plug-in rules for density level sets estimation}
\label{secDLSE}
\setcounter{equation}{0}
Fix $\lambda>0$ and recall that our goal is to estimate the connected
components $T_j=T_j(\lambda)$, $j=1,2, \ldots$, of $\Gamma= \Gamma(\lambda)=\{x \in \cX: p(x) \ge \lambda\}$, using
the unlabeled sample $\X_u$ of size $m$. A
simple and intuitive way to achieve this goal is to use \emph{plug-in
  estimators} of $\Gamma$ defined by
$$
\hat \Gamma= \hat \Gamma(\lambda)\triangleq  \left\{ x \in \cX: \hat p_m(x) \ge \lambda\right\}\,,
$$
where $\hat p_m$ is some estimator of $p$. A straightforward generalization
are the \emph{penalized plug-in estimators} of $\Gamma(\lambda)$, defined by
$$
\tilde \Gamma_{\ell}= \tilde \Gamma_{\ell} (\lambda) \triangleq \left\{ x \in
  \cX: \hat p_m(x) \ge \lambda+ \ell \right\}\,,
$$
where $\ell > 0$ is a penalization. Clearly $\tilde \Gamma_{\ell}
\subset \hat \Gamma$. Therefore the connected components of $\tilde \Gamma_\ell$ are farther from
each other than those of $\hat \Gamma$. Keeping in mind
that we want estimators that are consistent from inside we are going
to consider sufficiently large penalization $\ell=\ell(m)$.

Plug-in rules have a practical advantage over direct methods such as empirical
excess mass maximization (see, e.g., \cite{p95}, \cite{t97}, \cite{shs05}). Once
we have an estimator $\hat p_m$, we can compute the whole collection $\{\tilde
\Gamma_\ell(\lambda), \lambda>0 \}$, which might be of interest for the user who wants to try
several values of $\lambda$. Note also that a wide range of density
estimators is 
available in usual software. A density estimator can be parametric, typically
based on a mixture model,
or nonparametric such as histograms or kernel density estimators. 
\begin{defin}
For any $\lambda, \gamma \geq 0$,  a  function $f:\cX \to \R$ is said
to have $\gamma$-exponent at level $\lambda$ if there
exists a constant $c_0>0$ such that, for all $\varepsilon>0$,
$$
\Leb \left\{  |f(X)-\lambda| \leq \varepsilon\right\} \leq c_0 \varepsilon^\gamma\,.
$$
\end{defin}
It is an analog of the local margin assumption but for arbitrary level
$\lambda$ in place of $1/2$. When $\gamma>0$ it ensures that the function
$f$ has no flat part at level $\lambda$. 

The next theorem gives fast rates of convergence for penalized plug-in rules
when $\hat p_m$ satisfies an exponential inequality and $p$ has
$\gamma$-exponent at level $\lambda$. Moreover,
it ensures that when the penalization $\ell$ is suitably chosen, the plug-in estimator is consistent from inside.
\begin{TH1}
\label{fastrates}
Fix $\lambda>0, \gamma>0$ and $\Delta>0$. Let $\hat p_{m}$ be an estimator of the
density $p$ such that $P_X(\hat p_m(X) \ge
\lambda) \le C$, $\psubm$-almost surely for some positive constant $C$ and let
$\cP$ be a class of densities on $\cX$. Assume that there
exist positive constants $c_{1}, c_{2}$ and $a \le 1$, such that  for
$P_X$-almost all $x \in \cX$, we have
\begin{equation}
\label{expineq}
\sup_{p \in \cP} \psubm \left( |\hat p_{m}(x) - p(x)| \geq \delta \right) \leq
c_{1} \e^{-c_{2}m^a\delta^2}\,, \  m^{-a/2}< \delta< \Delta\,.
\end{equation}
Assume further that $p$ has $\gamma$-exponent at level $\lambda$  and  that the penalty $\ell$ is chosen as
\begin{equation}
\label{pen}
\ell= \ell(m)= m^{-\frac{a}{2}}\log m\,.
\end{equation}
Then the plug-in estimator $\tilde \Gamma_{\ell}$ is consistent from
inside at rate $m^{-\frac{\gamma a}{2}}$ . 
\end{TH1}
Consider a kernel density estimator $\hat p_m^K$ based on the
sample $\X_u$ defined by
\begin{equation}
\label{kde}
\hat p_{m}^K(x) \triangleq \frac{1}{mh^d} \sum_{i=n+1}^{n+m} K \left( \frac{X_{i}-
    x}{h} \right), \quad x \in \cX \,,
\end{equation}
where $h>0$ is the bandwidth parameter and $K: \cX \to \R$ is a kernel. If $p$
is assumed to have H\"older smoothness parameter $\beta>0$ and if  $K$ and $h$ are
suitably chosen, it is a standard exercise to prove inequality of type
\eqref{expineq} with
$a=2\beta/(2 \beta+d)$. In that case, it can
be shown that the rate $m^{-\frac{\gamma a}{2}}$ is optimal in a minimax sense.

\section{Discussion}
\label{secdisc}
We proposed a formulation of the cluster assumption in probabilistic terms. This formulation relies on Hartigan's
\cite{h75} definition of clusters but it can be modified to match other
definitions of clusters in the following way.
\begin{quote}
Consider a collection of $r_0$-connected and $s_0$-separated sets (clusters) $T_j, j=1,2,
\ldots$. Then the function
$x \mapsto (\eta(x)-1/2)$ has constant sign on each $T_j$.
\end{quote} 
We also proved that there is no hope to improve the classification performance
outside of these clusters. Based on these remarks, we defined the
$\lambda$-thresholded excess-risk which can be easily generalized to the setup
of general clusters defined above. 
Finally we proved that when we have consistent estimators of the
clusters, it is possible to achieve exponential rates of convergence for the
$\lambda$-thresholded excess-risk. The theory developed here can be extended to any
definition of clusters as long as they can be consistently estimated.

Note that our definition of clusters is parametrized by $\lambda$ which is left to the user, depending on his trust in the cluster
 assumption. The choice of $\lambda$ can be made by fixing
 $P_X(\Gamma^c)$, the probability of the rejection region. We refer to \cite{cff01}
 for more details. Note that data-driven choices of $\lambda$ could be
 easily derived if we impose a condition on the purity of the
 clusters, i.e. if we are given the $\delta$ in the global margin
 assumption. Such a choice could be made by decreasing $\lambda$ until
 the level of purity is attained. However, any data-driven choice of
 $\lambda$ has to be made using the labeled data. It would therefore
 yield much worse bounds. 

General open problems are: applying the cluster assumption to other definitions
of clusters and study the whole excess-risk in the framework of semi-supervised
classification with a rejection option. 

\section{Appendix: proofs}
\label{proofs}
\setcounter{equation}{0}
\subsection{Proof of Theorem~\ref{thpop}}
Using the decomposition of $\Gamma$ into its connected
components, we can decompose $\cE_\lambda(\hat g_n)$ into
$$
\cE_{\lambda}(\hat g_n)= \E_n\sum_{j\ge 1} \int_{T_j} |2 \eta(x) -1 | \1_{\{\hat g_{n}^j(x)\neq
  g^\star(x)\}}p(x) \ud x\,.
$$
Fix $j\in \{1,2, \ldots\}$ and assume w.l.o.g. that $\eta \ge 1/2$ on
$T_j$. It yields $g^\star(x)=1, \ \forall \, x \in T_j$, and since $\hat g_n$
is also constant on $T_j$,  we get
\begin{equation}
\label{prth11}
\begin{split}
\int_{T_j} |2 \eta(x) -1 | \1_{\{\hat g_{n}^j(x)\neq
  g^\star(x)\}}p(x) \ud x &=  \1_{\{Z_n^j \le 0\}}\int_{T_j} (2 \eta(x) -1
  )p(x) \ud x\\
& \le \delta_j\1_{\big\{|\delta_j-\frac{Z_n^j}{n}| \ge \delta_j\big\}}\,,
\end{split}
\end{equation}
Taking expectation $\E_n$ on both sides of \eqref{prth11} we get
\begin{equation}
\label{prth12}
\begin{split}
\E_n\int_{T_j} |2 \eta(x) -1 | \1_{\{\hat g_{n}^j(x)\neq
  g^\star(x)\}}p(x) \ud x& \le \delta_j
  \pn\Big[\big|\delta_j-\frac{Z_n^j}{n}\big| \ge \delta_j \Big]\\
& \le 2 \delta_j \e^{-n
  \delta_j^2/2}\,,
\end{split}
\end{equation}
where we used Hoeffding's inequality to get the last inequality. Summing now over $j$ yields the theorem.

\subsection{Proof of Proposition~\ref{propcc}}
\begin{figure}[h] \centering
\psfrag{r}{\Large $r_0$}
\psfrag{HH}{\Large $\tilde H_1$}
\psfrag{HT}{\Large $\tilde H_2$}
\psfrag{TJ}{\Large $T_j$}
\psfrag{lo}{\Large $(\log m)^{-1}$}
\begin{center}
\begin{minipage}[tr]{0.9\textwidth}
  \includegraphics[angle=0, width=0.7\textwidth]{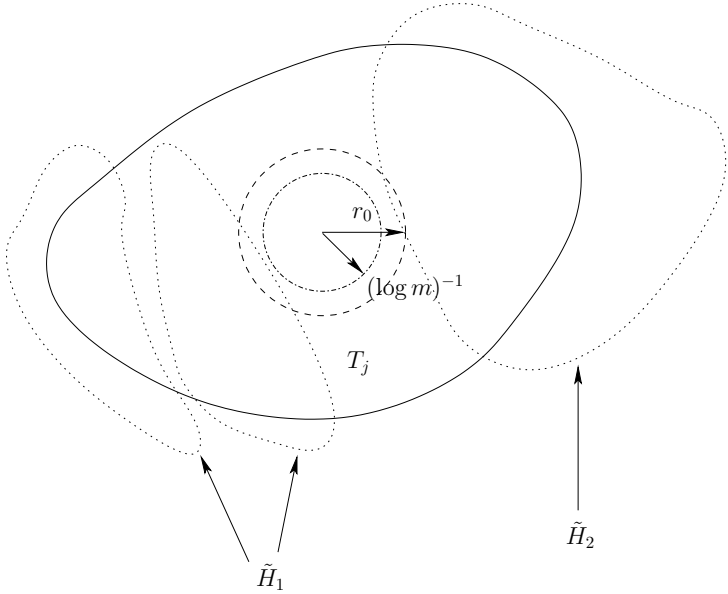}
\end{minipage}
\end{center}
\caption{By construction, $\tilde H_1$ and $\tilde H_2$ are separated by a ball of radius
  $(\log m)^{-1}$, which is included in $\cB(x,r_0)$ when $m \ge m_0$. So if $\{1,2\} \subset \kappa(j)$ or $\kappa(j)=\emptyset$, this ball is
  included in in $\tilde \Gamma_{\ell} \symdiff \Gamma $.}\label{kappa1}
\end{figure}

Define $m_0\triangleq
\exp(1/(r_0 \wedge s_0))$. Since the connected components $T_j$ are
$r_0$-connected, there is only a finite number $J\ge 1$ of them. We
simply denote $D(\cJ)$ by $D$. For any $j=1, \ldots, J$, the $r_0$
connectedness of $T_j$ yields on the one hand,
\begin{eqnarray*}
A_1(j) \triangleq \{{\rm card} [\kappa(j)]=0\} & \subset & \big\{\Leb \big[\tilde G_m \symdiff \Gamma\big]> \lambda c (\log
 m)^{-d} \big\}\,,\\
A_2(j) \triangleq\{{\rm card} [\kappa(j)]\ge 2\}& \subset &\big\{ \Leb\big[\tilde G_m \symdiff \Gamma\big]> \lambda c (\log
 m)^{-d} \big\}\,.
\end{eqnarray*}
The previous inclusions are illustrated in Figure~\ref{kappa1}.

\begin{figure}[h] \centering
\psfrag{r}{\Large $s_0$}
\psfrag{P}{\Large $\tilde T_l$}
\psfrag{S}{\Large $\tilde T_{l'}$}
\psfrag{T}{\Large $T_{j}$}
\psfrag{Q}{\Large $T_{j'}$}
\begin{minipage}[tr]{0.4\textwidth}
  \includegraphics[angle=0, width=\textwidth]{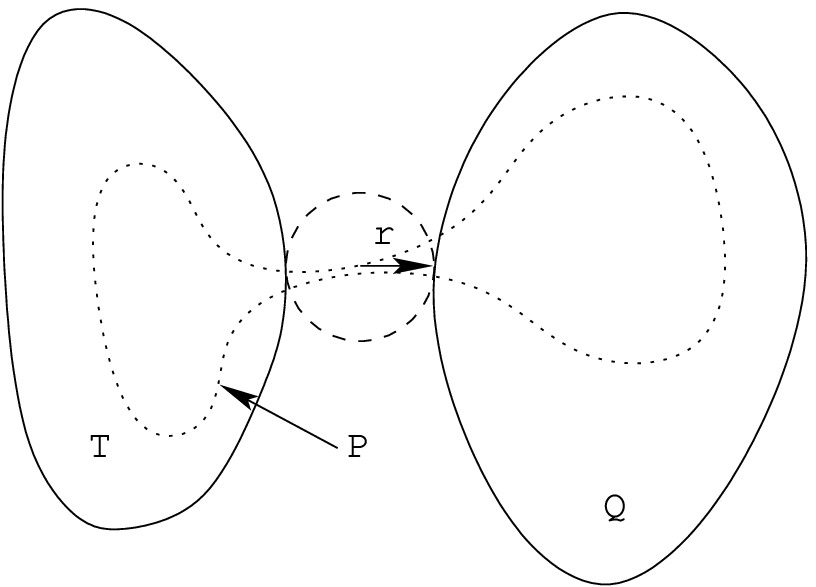}
\end{minipage}
\hspace{1cm}
\begin{minipage}[tr]{0.4\textwidth}
  \includegraphics[angle=0, width=\textwidth]{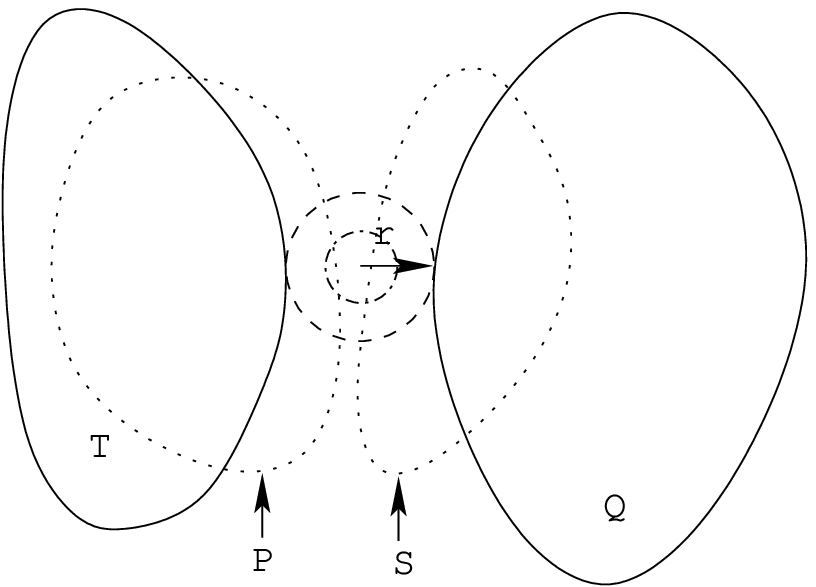}
\end{minipage}
\caption{Case $(i)$ (left) and case $(ii)$ (right).}\label{kappa2}
\end{figure}

\noindent On the other hand, $\kappa(j)\cap \kappa(j')\neq \emptyset$  for some $j' \neq
j$ when either $(i)$ $\exists \, l$ s.t. $\tilde
T_l \cap T_j \neq \emptyset, \ \tilde T_l \cap T_{j'} \neq \emptyset$ or
$(ii)$ $\exists \, l \neq l'$ s.t. $ \tilde
T_l \cap T_j \neq \emptyset, \ \tilde T_{l'} \cap T_{j'} \neq \emptyset$ and
$d_\infty(\tilde T_l, \tilde T_{l'}) < 2(\log m)^{-1}$. Both cases
yield the existence of $x\in \Gamma^c \cap \tilde G_m$ such that
$\cB(x,(\log m)^{-1}) \subset \Gamma^c$ for $m \ge m_0$. Therefore 
$$
\Leb(\tilde G_m \cap \Gamma^c) \ge \Leb(\tilde G_m \cap \cB(x,(\log m)^{-1}) )
$$
By construction of $\tilde
G_m$, we have $\Leb(\cB(x,(\log m)^{-1}) \cap \tilde G_m) \ge
 m^{-\alpha}(\log m)^{-d}$. Hence
$$
A_3(j) \triangleq\bigcup_{j' \neq j}\{\kappa(j) \cap\kappa(j') \neq
\emptyset\} \subset \big\{\Leb(\tilde G_m \cap \Gamma^c)\ge
 m^{-\alpha}(\log m)^{-d} \big\}
$$
Both cases are illustrated in Figure~\ref{kappa2}.

\noindent Now, since 
$$
D^c= \bigcup_{j=1}^J A_1(j) \cup A_2(j) \cup A_3(j),
$$
we get
$$
\psubm(D^c) \le \psubm \big\{\Leb \big[\tilde G_m \symdiff \Gamma\big]> \lambda c (\log
 m)^{-d} \big\}+\psubm \big\{\Leb(\tilde G_m \cap \Gamma^c)\ge
 m^{-\alpha}(\log m)^{-d} \big\}\,.
$$
Using the Markov inequality for both terms we obtain
$$
\psubm \big\{\Leb \big[\tilde G_m \symdiff \Gamma\big]>  \lambda c (\log
 m)^{-d} \big\}= \widetilde O \left(m^{-\alpha} \right)\,.
$$
and
$$
\psubm \big\{\Leb(\tilde G_m \cap \Gamma^c)\ge
m^{-\alpha}(\log m)^{-d} \big\}= \widetilde O \left(m^{-\alpha} \right)\,
$$
where we used the fact that $\tilde G_m$ is consistent from inside
with rate $m^{-\alpha}$. It yields the statement of the proposition.
\subsection{Proof of Theorem~\ref{mainth}}
The $\lambda$-thresholded excess-risk $\cE_\lambda(\tilde g_{n,m})$ can be
decomposed w.r.t the event $D$ and its complement. It yields
\begin{equation*}
\begin{split}
\cE_{\lambda}(\tilde g_{n,m})& \le \E_m\bigg[\1_{D}\E_n \Big(\int_{\Gamma}|2 \eta(x) -1 | \1_{\{\tilde g_{n,m}(x)\neq
  g^\star(x)\}}p(x) \ud x \Big| \X_u\Big)\bigg]+ \psubm\left(D^c \right)
\end{split}
\end{equation*}
We now treat the first term of the RHS of the above inequality, i.e., on the event $D$. Fix $j\in \{1,2, \ldots\}$ and assume w.l.o.g. that $\eta \ge 1/2$ on
$T_j$. Simply write $Z^k$ for $Z_{m,n}^k$. By definition of  $D$,
there is a one-to-one correspondence between the collection $\{T_j\}_j$ and
the collection $\{\tilde H_k\}_k$. We denote by $\tilde H_j$ the unique
element of $\{\tilde H_k\}_k$ such that $ \tilde H_j \cap T_j \neq
\emptyset$. On $D$, for any $j \ge 1$, we have,
\begin{equation*}
\begin{split}
&\E_n \Big( \int_{T_j} |2 \eta(x) -1 | \1_{\{\tilde g_{n,m}^j(x)\neq
  g^\star(x)\}}p(x) \ud x\Big| \X_u\Big)\\
&\phantom{mmmm} \le  \int_{T_j \setminus \tilde G_m} (2
\eta-1) \ud P_X+ \E_n \Big(\1_{\{Z^j \le 0\}}\int_{T_j \cap \tilde
  H_j}(2 \eta-1) \ud P_X\Big| \X_u\Big)\\
& \phantom{mmmm}\le L(p)\Leb (T_j \setminus \tilde G_m)+ \delta_j \pn \big(Z^j \le 0 |\X_u)
\end{split}
\end{equation*}
On the event $D$, For any $0< \theta < 1$, it holds
\begin{equation*}
\begin{split}
\pn\big(Z^{j} \le 0 |\X_u)  & = \pn\big(\int_{T_j}(2 \eta-1) \ud P_X-Z^{j}
\ge \delta_j |\X_u \big)\\
& \le \pn\big(\big|Z^{j} - \int_{\tilde H_{j}}(2 \eta-1) \ud P_X\big|
\ge \theta \delta_j  |\X_u \big)\\
&\phantom{mmmm} + \1_{\left\{P_X\big[T_j \symdiffsmall \tilde H_{j}\big] 
\ge (1-\theta)\delta_j \right\}}\,.
\end{split}
\end{equation*}
Using Hoeffding's inequality to control the first term, we get
$$
\pn\big(Z^{j} \le 0 |\X_u) \le 2 \e^{-n(\theta \delta_j)^2/2}+  \1_{\left\{P_X\big[T_j \symdiffsmall \tilde H_{j}\big] 
\ge (1-\theta)\delta_j \right\}}\,.
$$
Taking expectations, and summing over $j$, the $\lambda$-thresholded  excess-risk is upper bounded by
\begin{equation*}
\begin{split}
\cE_{\lambda}(\tilde g_{n,m})& \le  \frac{2L(p)}{1-\theta}\E_m\Big[ \Leb (\Gamma \symdiff \tilde G_m)\Big]
+ 2\sum_{j\ge 1} \delta_j \e^{-n (\theta \delta _j)^2/2}  + \psubm\left(D^c \right)\,,
\end{split}
\end{equation*}
where we used the fact that on $D$, 
$$
\sum_{j\ge 1} \Leb\big[T_j \symdiff \tilde H_{j}\big] \le \Leb\big[\Gamma
\symdiff \tilde G_m \big] \,.
$$
From Proposition~\ref{propcc}, we have $\psubm\left(D^c \right)=
\widetilde O \left(m^{-\alpha}
\right) $ and $\E_m\Big[ \Leb (\Gamma \symdiff \tilde G_m)\Big]=\widetilde O \left(m^{-\alpha}
\right)$ and the theorem is proved.

\subsection{Proof of Theorem~\ref{fastrates}}
Recall that
$$
\tilde \Gamma_{\ell} \symdiff \Gamma = \left( \tilde \Gamma_{\ell}
  \cap \Gamma^c \right) \sqcup \left( \tilde \Gamma_{\ell}^c \cap \Gamma
\right)\,.
$$
We begin by the first term. We have
$$
\tilde \Gamma_{\ell}
  \cap \Gamma^c=\big\{x \in \cX: \hat p_m(x) \ge \lambda + \ell, p(x)
 < \lambda  \big\}\subset \big\{x \in \cX: |\hat p_m(x) - p(x)| \ge \ell  \big\}\,.
$$
The Fubini theorem yields
$$
\E_m\big[\Leb (\tilde \Gamma_{\ell}
  \cap \Gamma^c ) \big] \le \Leb(\cX)\sup_{x \in \cX} \psubm \left[
  |\hat p_m(x) - p(x)| \ge \ell  \right] \le c_{3} \e^{-c_{2}m^a\ell^2}\,,
$$
where the last inequality is obtained using \eqref{expineq} and
$c_3=c_1 \Leb(\cX)>0$. Taking $\ell$ as in \eqref{pen} yields for $m
\ge \exp(\gamma a/c_2)$,
\begin{equation}
\label{prthfast1}
\E_m\big[\Leb (\tilde \Gamma_{\ell}
  \cap \Gamma^c ) \big] \le c_3m^{-\gamma a}.
\end{equation}
We now prove that 
$\E_m\big[\Leb (\tilde \Gamma_{\ell}
  \cap \Gamma^c ) \big] =\widetilde O \big(m^{-\frac{\gamma
    a}{2}}\big)$. Consider the following decomposition where we drop
the dependence in $x$ for notational convenience,
$$
 \tilde \Gamma_{\ell}^c \cap \Gamma=B_1 \cup B_2,
$$
where 
$$
B_1=\big\{ \hat p_m < \lambda + \ell, p
  \ge \lambda+2 \ell  \big\} \subset \big\{|\hat p_m - p| \ge \ell  \big\}
$$
and 
$$
B_2 = \big\{\hat p_m <
 \lambda + \ell, \lambda \le p(x)
< \lambda+2 \ell  \big\} \subset \big\{ |p-\lambda| \le \ell  \big\}\,.
$$
Using \eqref{expineq} and \eqref{pen} in the same fashion as above we
get $\E_m\big[\Leb (B_1 ) \big]=\widetilde O \big(m^{-\frac{\gamma
    a}{2}}\big)$. The term corresponding to $B_2$ is controlled using the
$\gamma$-exponent of density $p$ at level $\lambda$. Indeed, we have
$$
\Leb(B_2) \le c_0 \ell^\gamma \le c_0 (\log
m)^{\gamma}m^{-\frac{\gamma a}{2}}= \widetilde O \big( m^{-\frac{\gamma a}{2}}\big)
$$
The previous upper bounds for $B_1$ and $B_2$ together with
\eqref{prthfast1} yield the consistency from inside.

\end{document}